\documentclass[11pt]{amsart}
\usepackage{graphicx}
\usepackage{amssymb}
\usepackage{amsmath}
\usepackage{amsthm}
\usepackage{amscd}

\newtheorem{thm}{Theorem}[section]

\newtheorem{cor}[thm]{Corollary}
\newtheorem{lem}[thm]{Lemma}

\newtheorem{prop}[thm]{Proposition}
\theoremstyle{definition}

\theoremstyle{remark}
\newtheorem{rem}[thm]{Remark}
\numberwithin{equation}{section}

\begin{document}

\title[Duality between quantum symmetric algebras]{Duality between quantum symmetric algebras}
\author[Xiao-Wu Chen] {Xiao-Wu Chen }
\thanks{E-mail: xwchen$\symbol{64}$mail.ustc.edu.cn.}
\thanks{Supported by the National Natural Science Foundation of China (Grant No. 10301033 and No. 10501041).}
%%\subjclass{}%

%\keywords{Quantum groups, quivers}
%
\maketitle
%\date{}%
%\dedicatory{}%
%\commby{}%
% ----------------------------------------------------------------

\begin{center}
Department of Mathematics \\University of Science and Technology
of China \\Hefei 230026, Anhui, P. R. China
\end{center}

\vskip5pt

\begin{abstract}
Using certain pairings of couples, we obtain a large class of
two-sided non-degenerated graded Hopf pairings for quantum symmetric
algebras.
\end{abstract}

\section{Introduction}

One can construct three kinds of important graded Hopf algebras
$T_H(M)$, ${\rm Cot}_H(M)$ and $S_H(M)$ for a given Hopf algebra
$H$ and an $H$-Hopf bimodule $M$, which will be, respectively,
referred as the tensor Hopf algebra, the cotensor Hopf algebra and
the quantum symmetric algebra associated to the given couple $(H,
M)$. These constructions go back to Nichols [N], and they are
highlighted in Rosso's paper [Ro], who proves that the
non-negative part of the quantum enveloping algebra $U_q^{\geq
0}(\mathfrak{g})$ ($q$ not a root of unity) of a complex
semisimple Lie algebra $\mathfrak{g}$ is a quantum symmetric
algebra.\par \vskip 5pt

A special case of the above Hopf algebras is of particular interest
to representationists. If the Hopf algebra $H \simeq K \times \cdots
\times K$ as algebras, then $T_H(M)$ is a path algebra of some
quiver, and so the quantum symmetric algebra $S_H(M)$ is also a
quotient of the path algebra (e.g., see [C], [CR1] and [GS]).
Dually, if the Hopf algebra $H$ is group-like (not necessarily of
finite dimension), then ${\rm Cot}_H(M)$ is a path coalgebra [CM] of
some  quiver,  and hence $S_H(M)$ is a large subcoalgebra of the
path coalgebra. Note that both of the above observations will lead
to the concept of Hopf quivers by Cibils and Rosso [CR2]. In fact,
these quiver presentations of the  above Hopf algebras are very
useful to study certain Hopf algebras, see [CHYZ] and [OZ], and they
also could be used to classify some comodules of quantum groups, see
[CZ].\par \vskip 5pt

Inspired  by the above cited works, we study the three kinds of
graded Hopf algebras and their universal properties. Moreover, using
their universal properties, we can build up a large class of graded
Hopf pairings for the quantum symmetric algebras, which can be seen
as a generalization of a result by Nichols [N, Proposition 2.2.1]
and the self-duality of $U_q^{\geq 0}(\mathfrak{g})$ (and its
variants) via Rosso's isomorphism [Ro, Theorem 15] (and the remarks
thereafter).\par \vskip 5pt

The paper is organized as follows: section 2 is devoted to recall
the three known constructions of graded Hopf algebras from a given
couple. We also include their universal properties. In section 3, we
prove the main results of this paper: Theorem 3.1 and Theorem 3.2,
which claim that there exists graded Hopf pairings between certain
quantum symmetric algebras, and furthermore, the existence of
two-sided non-degenerated Hopf pairings characterizes quantum
symmetric algebras. Note that the proof uses the technical notion of
algebra-coalgebra pairings. As a special case, the self-duality of
Hopf algebras is briefly discussed in 3.7.\par \vskip 5pt

All algebras and coalgebras will be over a fixed field $K$, and
$\otimes $ means $\otimes_K$. Graded algebras (resp. coalgebras
$\cdots$) will always mean positively-graded algebras (resp.
coalgebras $\cdots$).

\section{Three constructions of graded Hopf algebras}
This section is devoted to fix some notation and to recall three
constructions of graded Hopf algebras from a given couple.

\subsection{}
Let us recall some basic definitions in Hopf algebras (see [Sw] and
[M]). Throughout $H$ will be a Hopf algebra with comultiplication
$\Delta_H$, counit $\varepsilon_H$ and antipode $S_H$. Sometimes we
denote the multiplication of $H$ by $m_H$ and the unit of $H$ by
$1_H$.
\par

An $H$-Hopf bimodule $M$ is an $H$-bimodule (with the $H$-actions
denoted by ``.'') and $H$-bicomodule with structure maps $\rho_l:
M \longrightarrow H \otimes M$ and $\rho_r: M \longrightarrow M
\otimes H$ such that
\begin{align*}
\rho_l(h.m.h')= \sum h_1m_{-1}h'_1 \otimes h_2.m_0.h'_2\\
\rho_r(h.m.h')= \sum h_1.m_{0}.h'_1 \otimes h_2m_1h'_2,
\end{align*}
where $h, h' \in H$ and $m \in M$, and we use the Sweedler notation,
i.e., $\Delta(h)=\sum h_1 \otimes h_2$, $\rho_l(m)=\sum
m_{-1}\otimes m_0$ and $\rho_r(m)=\sum m_0 \otimes m_{1}$ (e.g., see
[Sw, p.10 and p.32]). For example, $H$ itself  is an $H$-Hopf
bimodule with the regular bimodule structure and the bicomodule
structure maps $\rho_l=\rho_r=\Delta_H$.

 We will refer to the above
pair $(H, M)$ as a \textit{couple} throughout this paper, where $H$
is a Hopf algebra and $M$ an $H$-Hopf bimodule.
\par \vskip 5pt

\subsection{}A Hopf algebra  $H$ is said to be \textit{graded}, if there exists a decomposition of vector spaces
$H=\oplus_{n \geq 0} H_n$ such that
\begin{align*}
H_n H_m \subseteq H_{n+m},   \quad 1_H &\in H_0,\quad  S_H(H_n)  \subseteq H_n, \\
\Delta_H (H_n) \subseteq \sum_{i+j=n} H_i \otimes H_j, & \quad
\varepsilon_H(H_n)=0 \quad \mbox{if}\quad  n>0,
\end{align*}
where $n, m \geq 0$.
\par

Consider a  graded Hopf algebra H. It is clear that $H_0$ is a
subHopf algebra. Moreover, $H_1$ is an $H_0$-bimodule induced by
the multiplication inside $H$.  Note that
$$\Delta_H(H_1) \subseteq H_0 \otimes H_1 \oplus H_1 \otimes H_0, $$
thus there exist unique maps $\rho_l: H_1 \longrightarrow H_0 \otimes H_1$ and $\rho_r: H_1 \longrightarrow H_1 \otimes H_0$ such that
\begin{align*}
\Delta_H(m)= \rho_l(m)+ \rho_r(m), \quad \forall\;  m \in H_1.
\end{align*}
\par \vskip 5pt

The following result is well known and it can be easily checked.

\begin{lem}
Use the notation above. Then $H_1$ is an $H_0$-Hopf bimodule with
the structure maps $\rho_l$ and $\rho_r$, and thus $(H_0, H_1)$ is a
couple.
\end{lem}

We will say  that the above couple $(H_0, H_1)$ is
\textit{associated to} the graded Hopf algebra $H$.

%\begin{rem}
%Every component $H_n$ of $H$ is an $H_0$-Hopf bimodule. In fact,
%for $n >0$, first, $H_n$ is an $H_0$-bimodule, and  there exist
%unique maps $\rho_l: H_n \longrightarrow H_0 \otimes H_n$ and
%$\rho_r: H_n \longrightarrow H_n \otimes H_0$ such that
%\begin{align*}
%\Delta_H(m) \equiv \rho_l (m) + \rho_r(m) \quad (\mbox{ mod }
%\sum_{1 < i < n-1} {H_i \otimes H_{n-i}})
%\end{align*}
%for all $m \in H_n$. Thus $H_n$ becomes an $H_0$-Hopf bimodule.
%Certainly, $H_0$ itself is an $H_0$-Hopf bimodule.
%\end{rem}

\subsection{} We will
recall the constructions of graded Hopf algebras $T_H(M)$, ${\rm
Cot}_H(M)$ and $S_H(M)$ from a given couple $(H, M)$.

 \par \vskip 5pt

\subsubsection{$T_H(M)$} As in 2.1, $H$ is a Hopf algebra and $M$ an  $H$-Hopf bimodule.
Denote $T_H(M)$ the tensor algebra associated to the $H$-bimodule
$M$, i.e.,
\begin{align*}
T_H(M)= H \oplus M \oplus (M \otimes_H M) \oplus \cdots \oplus
M^{\otimes_H n} \oplus \cdots ,
\end{align*}
where $M^{\otimes_H n}= M^{\otimes_H {n-1}} \otimes_H M$ for
$n\geq 1$.\par \vskip 5pt

 To avoid confusion, we will write $T_H(M) \underline{\otimes}
T_H(M)$ for $T_H(M) \otimes T_H(M)$. Consider the following two
maps
$$\Delta_H: H \longrightarrow  H\underline{\otimes }H \quad \subseteq T_H(M) \underline{\otimes} T_H(M) \quad
 $$
 and
$$\rho_l \oplus \rho_r: M \longrightarrow H \underline{\otimes} M + M \underline{\otimes} H \quad \subseteq T_H(M) \underline{\otimes} T_H(M).$$
 It is direct to see that  $T_H(M) \underline{\otimes} T_H(M)$ is an $H$-bimodule via the algebra map $\Delta_H$
 and the map $\rho_l \oplus \rho_r$ is an $H$-bimodule morphism.
 Applying the universal property of the tensor algebra $T_H(M)$ (e.g., see [N, Proposition 1.4.1]),
 we obtain that there exists a unique algebra map
\begin{align*}
\delta: T_H(M) \longrightarrow T_H(M) \underline{\otimes} T_H(M)
\end{align*}
 such that $\delta|_H=\Delta_H$ and $\delta|_{M}= \rho_l \oplus
 \rho_r$.\par

  Using  a similar argument,
   we obtain a unique  algebra map
   \begin{align*}
   \epsilon: T_H(M) \longrightarrow K
   \end{align*}
    such that $\epsilon|_H=\varepsilon_H$ and $\epsilon|_M=0$.
    \par\vskip 5pt

 It is not hard to see that $(T_H(M), \delta, \epsilon)$ is a coalgebra,
 and thus $T_H(M)$ becomes a (graded) bialgebra. By [N, Proposition 1.5.1] or [M, Lemma 5.2.10], the bialgebra
$T_H(M)$ is a Hopf algebra. In fact, one can describe the antipode
more explicitly. Let $s_1: M \longrightarrow M$ be a map such that
 \begin{align*}
 s_1(m)= -\sum S_H(m_{
 -1}). m_0 . S_H(m_1),
 \end{align*}
  where  $({\rm Id}_H \otimes \rho_r) \rho_l(m)= \sum m_{-1} \otimes m_0 \otimes m_{1}$ and $m \in M$.
  Again by the universal property of the tensor algebra $T_H(M)$,
  there is a unique algebra map $s: T_H(M) \longrightarrow T_H(M)^{op}$
  such that $s|_H=S_H$ and $s|_M =s_1$, where $T_H(M)^{op}$ is the opposite algebra of
  $T_H(M)$. One can deduce that the map  $s$ is the antipode of $T_H(M)$
  (here, one may use [Sw, p.73, Ex. 2) ]).

    We will call the resulting Hopf algebra $T_H(M)$ the \textit{tensor Hopf algebra} associated to the couple $(H, M)$.
    \par\vskip 10pt

   We observe the following universal property of the tensor Hopf algebra $T_H(M)$.  The proof
   follows immediately from the universal property of the tensor algebras and then
   the coalgebra structure of $T_H(M)$.

    \begin{prop}
    Let $B=\oplus_{n \geq 0}{B_n}$ be a graded Hopf algebra with the associated couple $(B_0, B_1)$.
    Then there exists a unique graded Hopf algebra morphism $\pi_B: T_{B_0}(B_1) \longrightarrow B$
    such that the restriction of $\pi_B$ to $B_0\oplus B_1$ is the identity
    map.\par
Moreover, the map $\pi_B$ is surjective if and only if $B$ is
generated by $B_0$ and $B_1$.
    \end{prop}

%\noindent{\bf Proof}\quad By the universal property of the tensor
%algebra, we have a unique (graded) algebra morphism $\pi_B:
%T_{B_0}(B_1) \longrightarrow B$ such that $\pi_B|_{B_0\oplus B_1}$
%is the identity map. By [Sw, Lemma 4.0.4], it suffices to show
%that the map $\pi_B$ is a coalgebra morphism, i.e., for each $x
%\in T_{B_0}(B_1)$ we have
%\begin{align*}
%(\pi_B \otimes \pi_B)(\delta(x))= \Delta_B(\pi_B(x)).
%\end{align*}
%For this end, set $F=(\pi_B \otimes \pi_B)\circ \delta$ and $G=
%\Delta_B \circ \pi_B$. Since $\pi_B$ is an algebra morphism, thus
%both $F$ and $G$ are (graded) algebra morphisms from
%$T_{B_0}(B_1)$ to $B \otimes B$. Hence we need only to show $F=G$
%on a set of generators of the tensor algebra $T_{B_0}(B_1)$.
%However $F=G$ on the generating subset $B_0\oplus B_1$ trivially.
%This completes the proof. \hfill $\blacksquare$
%

\subsubsection{${\rm Cot}_H(M)$} This is the dual construction of 2.3.1. As before,  $(H, M)$ is a couple
    where $H$ is a Hopf algebra and $M$  an $H$-Hopf bimodule.\par

    Denote ${\rm Cot}_H(M)$ the cotensor coalgebra with respect to the $H$-bicomodule $M$
    (for details, see [N, p.1526], and [CHZ]), i.e.,
    \begin{align*}
    {\rm Cot}_H(M)= H \oplus M \oplus (M \square_H M) \oplus \cdots \oplus M^{\square_H n}\oplus \cdots,
    \end{align*}
    where $M \square_H M$ is the cotensor product of $M$, and
    we denote $M^{\square_H  n}= M \square_H M \square_H \cdots \square_H M$ (with $n$-copies of $M$).
     Note that $M^{\square_H  n}$ is a subspace of $M^{\otimes n}$,
      and elements $\sum m^{1} \otimes \cdots\otimes m^n \in M^{\otimes n}$
      which belong to $M^{\square_H n}$ will be written as $\sum m^1 \square_H \cdots \square_H
      m^n$.\par \vskip 5pt

    The coalgebra structure of ${\rm Cot}_H(M)$ is described as follows: the comultiplication
    $\Delta: {\rm Cot}_H(M) \longrightarrow {\rm Cot}_H(M) \otimes {\rm Cot}_H(M)$ is given
     by $\Delta|_H=\Delta$, $\Delta(m)= \rho_l(m)+ \rho_r(m)$ for all $m \in M$,
     and, in general,
   \begin{align*}
\Delta (\sum m^1 \square_H \cdots \square_H m^n)&=\sum (m^1)_{-1}
\otimes
( (m^1)_0 \square_H \cdots \square_H m^n)\\
&+ \sum_{i=1}^{n-1} (m^1 \square_H \cdots \square_H m^i) \otimes
(m^{i+1} \square_H  \cdots \square m^n)\\
&+ \sum( m^1 \square \cdots \square_H (m^n)_0) \otimes (m^n)_1 \\ &
\in (H \otimes M^{\square_H n})\\ & \bigoplus_{i=1}^{n-1}(M^{\square_H
i} \otimes M^{\square_H
  (n-i)}) \oplus( M^{\square_H n }\otimes
 H)
 \\
& \subseteq {\rm Cot}_H{M} \otimes {\rm Cot}_H(M),
\end{align*}
for any $\sum m^1 \square_H \cdots \square_H m^n\in M^{\square n}$.
The counit $\varepsilon: {\rm Cot}_H(M) \longrightarrow K$ is given
by $\varepsilon |_H=\varepsilon_H$ and $\varepsilon |_{M^{\square_H
n}}=0$ if $n \geq 1$.\par \vskip 5pt

Endow a bialgebra structure on ${\rm Cot}_H(M)$ as follows: by the
universal property of the cotensor coalgebra ${\rm Cot}_H(M)$
(e.g., see [N, Proposition 1.4.2] or  [CHZ, Lemma 3.2]), there
exist unique graded coalgebra morphisms
\begin{align*}
m : {\rm Cot}_H(M) \otimes {\rm Cot}_H(M) \longrightarrow {\rm Cot}_H(M) \quad \mbox{and} \quad e: K \longrightarrow {\rm Cot}_H(M),
\end{align*}
such that $m|_{H \otimes H}$ is just the multiplication $m_H$ of
the Hopf algebra $H$ and $m|_{(H\otimes M) \oplus (M\otimes H)}$is
given by
\begin{align*}
m(h\otimes n + n' \otimes h')=h.n + n'.h'
\end{align*}
for all $n, n' \in M$ and $h, h' \in H$, and the unit map $e: K
\longrightarrow {\rm Cot}_H(M)$ maps  $1_K$ to $1_H$.
\par \vskip 10pt

One can verify  that  $({\rm Cot}_H(M), m, e)$ is an algebra, and
thus ${\rm Cot}_H(M)$ becomes a (graded) bialgebra. By [N,
Proposition 1.5.1] or [M, Lemma 5.2.10], the bialgebra ${\rm
Cot}_H(M)$ is a graded Hopf algebra. Denote its antipode by $S$. It
is not hard to see that
 \begin{align*}
 S|_H=S_H\quad \mbox{and}\quad  S(m)=
-\sum S_H(m_{
 -1}). m_0 . S_H(m_1)
 \end{align*}
  for all $m \in M$. In fact, by the universal property of the cotensor
coalgebra ${\rm Cot}_H(M)$, the graded anti-coalgebra morphism $S$
is uniquely determined by the above two identities.

The resulting Hopf algebra $({\rm Cot}_H(M), m, e, \Delta,
\varepsilon, S)$ will be
  called the \textit{cotensor Hopf algebra} associated to the couple $(H, M)$.\par \vskip 5pt

Recall that in a coalgebra $(C, \Delta_C, \varepsilon_C)$, the
wedge is defined as $V \wedge_C W= \Delta_C^{-1}(V \otimes C + C
\otimes W)$ for any subspaces $V$ and $W$ of $C$. An important
fact is that $H \wedge_{{\rm Cot}_H(M)} H=H\oplus M$. (To see
this, first note that $H\oplus M \subseteq H \wedge_{{\rm
Cot}_H(M)} H$. Since ${\rm Cot}_H(M)$ is a graded coalgebra, then
$H \wedge_{{\rm Cot}_H(M)} H$ is a graded subspace. Therefore, it
suffices to show that
 $$(H \wedge_{{\rm Cot}_H(M)} H)\cap
M^{\square_H n}=0, \quad n \geq 2.$$ In fact, let $x\in (H
\wedge_{{\rm Cot}_H(M)} H)\cap M^{\square_H n} $ be a nonzero
element. Thus $\Delta(x)\in H \otimes M^{\square_H n} + M^{\square_H
n} \otimes H$. Note that $x \in M^{\square_H n}$ and by the
definition of the comultiplication $\Delta$, we know that the term
belonging to $M \otimes M^{\square_H {n-1}}$ which occurs in
$\Delta(x)$ is not zero. This is a contradiction.)\par \vskip 5pt

 Dual to Proposition 2.2, we observe the following universal property of the cotensor Hopf algebra ${\rm Cot}_H(M)$.
  Note that it is a slight generalization of [OZ, Theorem 4.5] (the proof of the second statement
  needs to use the above recalled fact).

 \begin{prop}
 Let $B=\oplus_{n \geq 0}{B_n}$ be a graded Hopf algebra with the associated couple $(B_0, B_1)$.
 Then there exists a unique graded Hopf algebra morphism $i_B: B \longrightarrow {\rm Cot}_{B_0}(B_1) $
 such that its restriction  to $B_0\oplus B_1$ is the identity
 map.\par
 Moreover, the map $i_B$ is injective if and only if $B_0 \wedge_B
 B_0=B_0\oplus B_1$.
 \end{prop}

 %\noindent{\bf Proof.} \quad First by  the universal property of the cotensor coalgebra
% ${\rm Cot}_H(M)$ (e.g., see [N, Proposition 1.4.2] or [CHZ, Lemma 3.2])
% there exists a unique graded coalgebra morphism $i_B: B \longrightarrow {\rm
% Cot}_{B_0}(B_1)$ such that its restriction to $B_0\oplus B_1$ is the identity
% map. By [Sw, Lemma 4.0.4], it suffices to show that $i_B$ is an
% algebra morphism.\par
% Indeed, we have to show that $m\circ (i_B \otimes i_B)=
% i_B \circ
% m_B$, where $m$ and $m_B$ denote the multiplication maps of ${\rm
% Cot}_H(M)$ and $B$, respectively.  To see this, set $F= m\circ (i_B \otimes i_B)
% $ and $G= i_B \circ m_B$. Since $i_B$ is a coalgebra morphism, we
% deduce that both $F$ and $G$ are coalgebra morphisms from $B \otimes
% B$ to ${\rm Cot}_H(M)$. Note that $\pi_i\circ F=\pi_i \circ G$,
% $i=0, 1$, where $\pi_0$ and $\pi_1$ denote the projection maps from ${\rm Cot}_H(M)$
%to $H$ and $M$, respectively. Now apply the universal property of
%the cotensor coalgebra ${\rm Cot}_H(M)$ again, we deduce that
%$F=G$. This proves the first statement.\par
%
% For the second, first the ``only if '' part follows from the fact
% $B_0
%\wedge_{{\rm Cot}_{B_0}(B_1)} B_0=B_0\oplus B_1$ directly. To see
%the "if" part, denote by $\{B_{(i)}\}_{i\geq 0}$ the coradical
%filtration of $B$. By [Sw, Proposition 11.1.1], we have $B_{(0)}
%\subseteq B_0 $, and thus by the assumption
%$B_{(1)}=B_{(0)}\wedge_B B_{(0)} \subseteq B_0 \wedge_B
%B_0=B_0\oplus B_1$. Now the fact that $i_B$ is injective follows
%from Heynemann-Radford Lemma ([M, Theorem
% 5.3.1]). \hfill $\blacksquare$

 \subsubsection{$S_H(M)$} As above, $(H, M)$ is a couple. Denote $S_H(M)$ be the graded
 subalgebra of the cotensor Hopf algebra ${\rm Cot}_H(M)$ generated by $H$ and $M$. Clearly, $S_H(M)$ is
 a graded subHopf algebra of ${\rm Cot}_H(M)$.
 We will call $S_H(M)$ the \textit{quantum symmetric algebra} (see [Ro, p.407])
 associated to the couple $(H, M)$. \par \vskip 5pt

The following result is a direct consequence of Proposition 2.3.

 \begin{cor}
 Let $B= \oplus_{n \geq 0} B_n$ be a graded Hopf algebra generated by $B_0$ and $B_1$.
 Denote $(B_0, B_1)$ the couple associated to $B$.
 Then there exists a unique graded Hopf algebra epimorphism $j_B: B \longrightarrow S_{B_0}(B_1)$
 such that its restriction to ${B_0 \oplus B_1}$ is the identity map.
 \end{cor}
%
%\noindent {\bf Proof.}\quad By Theorem 2.5, we have a graded Hopf
%algebra morphism $i_B: B \longrightarrow {\rm Cot}_{B_0}(B_1)$
%such that its restriction to $B_0 \oplus B_1$ is the identity map.
%Since $B$ is generated by $B_0$ and $B_1$, it follows that the
%image of $i_B$ is exactly $S_{B_0}(B_1)$. Thus the map $i_B$
%induces the required epimorphism $j_B$. \hfill $\blacksquare$

 \begin{rem}
 (1).\quad We may deduce Nichols's result [N, 2.2] from Proposition 2.3:  the map  $i_B$ in Proposition 2.3
  is an isomorphism
 if and only if  $B$ is generated by $B_0$ and $B_1$ and $B_0 \wedge_B B_0= B_0\oplus B_1$.
 In particular, if $B_0$ is cosemisimple, we see that $i_B: B \simeq S_{B_0}(B_1)$ if and only
 if $B$ is coradically-graded and generated by $B_0$ and $B_1$. \par
 (2).\quad As a special case of Corollary 2.4, for every couple $(H, M)$,
 there is a unique graded Hopf algebra
  epimorphism
  \begin{align*}
  \pi_{(H, M)}: T_H(M) \longrightarrow S_H(M)
  \end{align*}
  such that its restriction to ${H \oplus M}$ is the identity map, where we denote $j_{T_H(M)}$ by
  $\pi_{(H, M)}$. Denote the kernel of $\pi_{(H, M)}$ by $I(H, M)$, hence it is a graded Hopf ideal of $T_H(M)$ and $T_H(M)/{I(H, M)} \simeq S_H(M)$.
 \end{rem}

 \section{Duality between quantum symmetric algebras}

 \subsection{} Let us recall the definition of Hopf pairings (see [K], p.110).
  Let $H$ and $B$ be Hopf algebras.
 A Hopf pairing $\phi: H \times B \longrightarrow K$ is a bilinear map such that
 \begin{align*}
 &\phi(1_H, b)=\varepsilon_B (b); \quad  \phi(h,1_B)=\varepsilon_H(h);\\
 &\phi(h, bc)=\sum \phi(h_{1}, b) \phi(h_{2}, c);\\
 &\phi(hg, b)=\sum \phi(h, b_{1}) \phi(g, b_{2});\\
 &\phi(S_H(h), b)=\phi(h, S_B(b)),
 \end{align*}
 where $h, g \in H$, $b, c \in B$, and $S_H$ and $S_B$ are the antipodes of $H$ and
 $B$, respectively. Note that one can define the \textit{transpose} $\phi^t: B \times H \longrightarrow K$ by $\phi^t(b, h)=\phi(h, b)$,
 which is also a Hopf pairing.\par \vskip 5pt

 Assume further that $H=\oplus_{n \geq 0} H_n$ and $B=\oplus_{n \geq 0}B_n$ are graded Hopf algebras.
  A Hopf pairing $\phi: H \times B \longrightarrow K$ is said to be
  \textit{graded}
  if $\phi(H_i, B_j)=0$ if $i \neq j$.\par \vskip 5pt

 \subsection{} We introduce an analogous concept of Hopf pairings. Let $(H, M)$ and $(B, N)$ be
 couples. Denote the $H$-comodule ({\it resp.} $B$-comodule) structure
 on $M$ ({\it resp.} $N$) by $\rho_l$ and $\rho_r$ ({\it resp.} $\delta_l$ and
 $\delta_r$). Denote the actions by ``.'' . \par \vskip 5pt

  A \textit{pairing} between the couples $(H, M)$ and $(B, N)$, denoted by
  \begin{align*}
  (\phi_0, \phi_1): (H, M)\times (B, N) \longrightarrow K,
  \end{align*}
   is given by a Hopf pairing $\phi_0: H \times B \longrightarrow K$ and
   a bilinear map $\phi_1: M \times N \longrightarrow K$ such that
 \begin{align*}
  \phi_1(h.m.g, n)&=\phi_0(h, n_{-1}) \phi_1(m, n_{0}) \phi_0(g, n_{1}),\\
  \phi_1(m, b.n.c)&=\phi_0(m_{-1}, b)\phi_1(m_{0}, n)\phi_0(m_{1}, c),
  \end{align*}
  where $h, g \in H$, $b, c\in B$, $m \in M$ and $n \in N$,
   and $({\rm Id}_{H}\otimes \rho_r) \rho_l(m)=\sum m_{-1} \otimes m_0 \otimes m_1 $
   and $({\rm Id}_B \otimes \delta_r) \delta_l(n)=\sum n_{-1} \otimes n_0 \otimes n_1 $.
   \par \vskip 10pt

 We have our main results.
  \begin{thm}
 Let $(\phi_0, \phi_1): (H, M)\times (B, N) \longrightarrow K$ be a pairing between couples.
  Then there exists a unique graded Hopf pairing
  \begin{align*}
  \phi: S_H(M) \times S_B(N) \longrightarrow K
  \end{align*}
  extending $\phi_0$ and $\phi_1$.\par \vskip 3pt
   Moreover, $\phi$ is two-sided non-degenerated
  if and only if $\phi_0$ and $\phi_1$ are.
\end{thm}

and

\begin{thm}
Let $H=\oplus_{n \geq 0} H_n$(resp. $B=\oplus_{n \geq 0}B_n$) be
graded Hopf algebras generated by $H_0$ and $H_1$ (resp. $B_0$ and
$B_1$). Assume that there exists a  two-sided non-degenerated graded
Hopf pairing $\psi: H \times B \longrightarrow K$. Then $j_H:H
\simeq S_{H_0}(H_1)$ and $j_B: B \simeq S_{B_0}(B_1)$, where the
maps $j_H$ and $j_B$ are explained in Corollary 2.4.
\end{thm}

\vskip 10pt

  \subsection{}To prove the above two results, we need to introduce the following technical concept,
   which is essentially the same as the (graded) duality between algebras and coalgebras.\par \vskip 5pt

  Let $A$ be an algebra and $(C, \Delta_C, \varepsilon_C)$ a
  coalgebra.
 Let $\phi: A \times C \longrightarrow K$ be a bilinear map,
  and
  define $\phi^*: A \longrightarrow C^*$ by $\phi^*(a)(c)=\phi(a, c)$.
  We say that $\phi$ is an \textit{algebra-coagebra pairing} if
  \begin{align*}
  \phi(1_A, c)=\varepsilon_C(c) \quad \mbox{and} \quad \phi(aa', c)=\sum \phi(a, c_1) \phi(a, c_2),
  \end{align*}
  for all $a, a' \in A$ and $c \in C$, where $1_A \in A$ is the identity element and
   $\Delta(c)=\sum c_1 \otimes
  c_2$.\par \vskip 5pt

  In fact, it is easily checked that $\phi$ is an algebra-coalgebra pairing
  if and only if $\phi^*$ is an algebra morphism, where $C^*$ is the dual algebra of the coalgebra $C$.
  \par \vskip 5pt

  The graded version of the above concept is as follows:
  let $A=\oplus_{n \geq 0}A_n$ be a graded algebra and $C=\oplus_{n \geq 0}C_n$ a graded coalgebra,
   an algebra-coalgebra pairing $\phi: A \times C \longrightarrow K$ is said to be \textit{graded},
   if $\phi(A_i, C_j)=0$ for $i \neq j$. As above, we can
   define a graded map $\phi^*: A \longrightarrow C^{\rm gr}$,
   where $C^{\rm gr}=\oplus_{n \geq 0} C_n^*$ is the graded dual of $C$.
  One sees that $\phi$ is a graded algebra-coalgebra pairing
  if and only if $\phi^*$ is a graded algebra map.\par \vskip 10pt

  In what follows, we assume that $A$ is an algebra and $M$ an $A$-bimodule (with actions denoted by ``.''),
  and $(C=\oplus_{n \geq 0}C_n, \Delta_C, \varepsilon_C)$ is a graded coalgebra.
  Clearly there exist unique maps
  \begin{align*}
  \delta_l: C_1 \longrightarrow C_0 \otimes C_1 \quad \mbox{and} \quad \delta_r:C_1 \longrightarrow C_1 \otimes C_0
  \end{align*}
  such that $\Delta_C(c)=\delta_l(c) + \delta_r(c)$ for all $c\in C_1$. By abuse of notation,
  write $\delta_l(c)= \sum c_{-1}\otimes c_{0}$ and $\delta_r(c)=\sum c_{0} \otimes c_1$. Thus $\Delta_C(c)=\sum c_{-1} \otimes c_0 +
  \sum c_0 \otimes c_1$.\par \vskip 5pt

  Assume further that $\phi_0: A \times C_0 \longrightarrow K $ is an algebra-coalgebra pairing,
  and $\phi_1: M \times C_1 \longrightarrow K$ is a bilinear map such that
  \begin{align}
  \phi_1(a.m.a', c)=\sum \phi_0(a, c_{-1}) \phi_1(m, c_0) \phi_0(a', c_1)
  \end{align}
  for all $a, a' \in A$ and $c \in C_1$, where $({\rm Id}_{C_0} \otimes \delta_r)\delta_l(c)=\sum c_{-1} \otimes c_0 \otimes c_1$.
  \par \vskip 10pt

  We have the following

  \begin{lem}
  Assume that $\phi_0$ and $\phi_1$ are as above. There exists a unique graded algebra-coalgebra pairing
   $\phi: T_A(M) \times C \longrightarrow K$ extending $\phi_0$ and $\phi_1$.
  \end{lem}

  \noindent {\bf Proof}\quad This is just a variant of the universal property of the tensor algebra $T_A(M)$.
   Using $\phi_0$ and $\phi_1$, we can define $\phi_0^*: A \longrightarrow C_0^* \subseteq C^{\rm gr}$ and
   $\phi_1^*: M \longrightarrow C_1^* \subseteq C^{\rm gr}$. Note that $\phi_0^*$ is an algebra map and
   $\phi_1^*$ is an $A$-bimodule morphism (exactly by the condition (3.1)).
   \par \vskip 3pt
  Now by the universal property of the tensor algebra $T_A(M)$, there exists a unique graded algebra map
   \begin{align*}
  \phi^*: T_A(M) \longrightarrow C^{\rm gr}
  \end{align*}
  extending $\phi_0^*$ and $\phi_1^*$. Define $\phi$ by $\phi(x, c)=\phi^*(x)(c)$, for all $x\in T_A(M)$ and $c \in C$.
  Immediately, $\phi$ is the unique graded algebra-coalgebra pairing extending $\phi_0$ and $\phi_1$. This completes the proof. \hfill $\blacksquare$
  \vskip 10pt

  \subsection{} Recall that any pairing
  $\phi: A \times C \longrightarrow K$ is said to be \textit{left non-degenerated }provided that
  for each nonzero $y\in C$ there is some $x \in A$ such that $\phi(x ,y) \neq 0$. Let us go to the situation of Theorem 3.1 and 3.2:
  we are given a pairing of couples $(\phi_0, \phi_1): (H, M)\times (B, N) \longrightarrow K$.
  The following result is of independent interest. \par \vskip 10pt

  \begin{prop}
  There exists a unique graded Hopf pairing
  \begin{align*}\phi: T_H(M) \times {\rm Cot}_B(N)\longrightarrow K
  \end{align*}
   extending $\phi_0$ and $\phi_1$. \par \vskip 3pt
   Moreover, if $\phi_0$ and $\phi_1$ are left non-degenerated, then so is $\phi$.
  \end{prop}

  \noindent{\bf Proof}\quad By Lemma 3.3, there exists a unique graded algebra-coalgebra pairing
  $\phi: T_H(M) \times {\rm Cot}_B(N)\longrightarrow K$ extending $\phi_0$ and $\phi_1$.
  We will  show that $\phi$ is a Hopf pairing. \par \vskip 5pt

  Use the notation in 2.3.1 and 2.3.2. First
   we have $\phi(x,1_B)=\epsilon(x)$ and
   $\phi(1_H, y)=\varepsilon(y)$ for all $x \in T_H(M)$, $y\in {\rm
   Cot}_B(N)$. (To see this, since $\phi$ is graded, we have $\phi(x, 1_B)=0=\epsilon(x)$
   for $x \in M^{\otimes _H n}$,
   $n \geq 1$; and for $x \in H$, $\phi(x, 1_B)=\phi_0(x, 1_B)=\epsilon(x)$.
    Similarly one obtains that $\phi(1_H, y)=\varepsilon(y)$.)
   Define two bilinear maps
  \begin{align*}
\Psi, \Phi: T_H(M) \times ({\rm Cot}_B(N)\otimes {\rm Cot}_B(N))
\longrightarrow K
  \end{align*}
  such that $\Psi(x, y \otimes z)= \phi(x, yz)$ and $\Phi(x, y \otimes z)=\sum \phi(x_1, y)\phi(x_2, z)$,
   $x\in T_H(M)$, $y, z\in {\rm Cot}_B(N)$. Note that both $\Psi$ and $\Phi$ are
   graded algebra-coalgebra pairings,
   and by the defining properties of the pairing $(\phi_0, \phi_1)$, we have
  \begin{align*}
  \Psi|_{H \times (B\otimes B)}=\Phi|_{H \times (B\otimes B)} \quad \mbox{and} \quad \Psi|_{M \times (B\otimes N+ N \otimes B)}=\Phi|_{M \times (B\otimes N+ N \otimes B)}.
  \end{align*}
  Now by the uniqueness part of Lemma 3.3,
   we obtain that $\Psi=\Phi$, i.e., $\phi(x, yz)=\phi(x_1, y)\phi(x_2,
   z)$.\par \vskip 3pt

  Similarly, we construct two graded algebra-coalgebra pairings
  \begin{align*}
  \Psi', \Phi': T_H(M) \times {\rm Cot}_B(N)^{cop}\longrightarrow K
  \end{align*}
  such that $\Psi'(x, y)=\phi(s(x), y)$ and $\Phi'(x, y)=\phi(x, S(y))$,
  where $x\in T_H(M)$ and $y \in {\rm Cot}_B(N)$, and ${\rm Cot}_B(N)^{cop}$ denotes  the  opposite coalgebra.
  By a similar argument as above, we show that $\phi(s(x), y)=\phi(x, S(y))$.
  \par \vskip 3pt

  Summing up the above, we have shown that $\phi$ is the unique required
   graded Hopf pairing. \par \vskip 5pt

  For the second statement, assume that $\phi_0$ and $\phi_1$ are left
  non-degenerated,
   we need to show that for every nonzero element $y \in N^{\square_B i}$,
   there exists some $x \in M^{\otimes_H i}$ such that $\phi(x, y)\neq 0$, $i \geq 2$.
   Since   $\phi_1: M \times N \longrightarrow K$ is left non-degenerated,
   hence the following bilinear map will be left non-degenerated:
  \begin{align*}
 \phi_1^{\otimes i}: M^{\otimes i} \times N^{\otimes i} \longrightarrow K,
  \end{align*}
 where $\phi_1^{\otimes i}(m^1 \otimes \cdots \otimes m^i, n^1 \otimes \cdots \otimes n^i )=
 \prod_{r=1}^i \phi_1(m^r, n^r)$. Note that $N^{\square_B i} \subseteq N^{\otimes i}$,
 hence for the nonzero $y \in N^{\square_B i}$,
 there exists some $x' \in M^{\otimes i}$ such that $\phi_1^{\otimes i} (x', y)\neq 0$.
 \par \vskip 3pt

 Denote by  $p: M^{\otimes i} \longrightarrow M^{\otimes_H i}$ the natural projection
 map. By the fact that $\phi$ is an algebra-coalgebra pairing, we have
  \begin{align*}
  \phi(p(x'), y)=\phi_1^{\otimes n}(x', y).
  \end{align*}
 Take  $x=p(x')$. We see that $\phi(x, y) \neq 0$, finishing the proof. \hfill $\blacksquare$
 \par \vskip 15pt

 \subsection{Proof of Theorem 3.1:}\quad
 Consider the following composite of   morphisms between graded
 Hopf algebras
 \begin{align*}
 \pi: T_B(N) \stackrel{\pi_{(B, N)}}{\longrightarrow} S_B(N) \hookrightarrow {\rm Cot}_B(N),
 \end{align*}
 where the map $\pi_{(B, N)}$ is described in Remark 2.5(2) and the second map is just the inclusion.
  Applying Proposition 3.4, we have a graded Hopf pairing $\phi': T_H(M) \times {\rm Cot}_B(N) \longrightarrow
  K$ extending $\phi_0$ and $\phi_1$. Define
 \begin{align*}
 \phi'': T_H(M) \times T_B(N) \longrightarrow K.
 \end{align*}
by putting $\phi''(x, y)= \phi'(x, \pi(y))$. Thus $\phi''$ is a
graded Hopf pairing.
\par \vskip 5pt

 Note that $I(B, N)$ is the kernel of  $\pi_{(B,N)}$ and thus the kernel of $\pi$,
 we see that
 \begin{align*}
 \phi''(T_H(M), I(B, N))=0.
 \end{align*}
 \par \vskip 5pt

  We now claim that
  \begin{align*}
  \phi''(I(H, M), T_B(N))=0.
  \end{align*}

For this end, apply Proposition 3.4 again, we have a graded Hopf
pairing
 \begin{align*}
 \psi: {\rm Cot}_H(M) \times T_B(N) \longrightarrow K
 \end{align*}
 extending $\phi_0$ and $\phi_1$.  Consider the following composite
  \begin{align*}
 \pi': T_H(M) \stackrel{\pi_{(H, M)}}{\longrightarrow} S_H(M) \hookrightarrow {\rm
 Cot}_H(M).
 \end{align*}
 Define $\psi': T_H(M) \times T_B(N) \longrightarrow K$ by $ \psi'(x, y)=\psi(\pi'(x),
 y)$. Since $\pi'$ is a (graded) Hopf algebra morphism, thus
 $\psi'$ is a graded Hopf pairing. Similarly as above, we have $ \psi'(I(H, M), T_B(N))=0.$
 Note that both $\phi''$ and $\psi'$ are graded algebra-coalgebra
 pairings extending $\phi_0$ and $\phi_1$. Applying Lemma 3.3, we
 have $\phi''=\psi'$. This proves the claim.\par

  \par \vskip 5pt
  So we have shown that $\phi''(T_H(M), I(B, N))=0$ and $\phi''(I(H, M), T_B(N))=0$. Recall from
  Remark 2.5(2) that we have
  \begin{align*}
  T_H(M)/{I(H, M)} \simeq S_H(M)\quad \mbox{and}\quad  T_B(N)/{I(B, N)} \simeq
  S_B(N).
  \end{align*}
  Thus we deduce that
  $\phi''$ induces a unique graded Hopf pairing
  $$\phi: S_H(M) \times S_B(N) \longrightarrow
  K$$
   such that the following diagram commutes
 \[
\begin{CD}
T_H(M) \times T_B(N) @>  \phi''>> K    \\
@V \pi_{(H, M)}\times \pi_{(B, N)} VV   @V {\rm Id}_K VV      \\
S_{H}(M) \times S_{B}{(N)} @> \phi >>K.
\end{CD}\]

  Explicitly, $\phi(\pi_{(H, M)}(x), \pi_{(B, N)}(y))=\phi''(x, y)$, for all $x \in T_H(M)$ and $y \in
  T_B(N)$.\par \vskip 3pt
  Obviously, the pairing $\phi$ extends the maps $\phi_0$ and $\phi_1$, as required.
  Note that the uniqueness of $\phi$ is trivial,
  since $S_H(M)$, as an algebra, is generated by $H$ and $M$. (Here, one needs to consult the fourth identity
  in the definition of Hopf pairing, see 3.1). \par \vskip 5pt

  For the second statement, assume that $\phi_0$ and $\phi_1$ are two-sided non-degenerated.
  By Proposition 3.4, we have that $\phi'$ is left non-degenerated.
  Note that
  $\phi(\pi_{(H, M)}(x), \pi_{(B, N)}(y))=\phi''(x, y)=\phi'(x, \pi(y))$.
  This implies that $\phi$ is left non-degenerated. For right non-degeneratedness, first apply Proposition 3.4 to $\psi^t$
  (the transpose of $\psi$, see 3.1),
   we deduce that $\psi^t$ is left  non-degenerated, that is, $\psi$ is right non-degenerated. Now
   note that
  $\phi(\pi_{(H, M)}(x), \pi_{(B, N)}(y))=\phi''(x, y)=\psi'(x, y)=\psi(\pi'(x), y)$,
  which implies that $\phi$ is right non-degenerated. This completes the proof. \hfill $\blacksquare$

  \vskip 15pt

  \subsection{Proof of Theorem 3.2:}\quad Since the Hopf pairing $\psi: H \times B\longrightarrow K$ is two-sided
  non-degenerated,
  so are the restrictions $\phi_0:=\psi|_{H_0 \times B_0}$ and $\phi_1:=\psi|_{H_1 \times B_1}$.
  Now applying  Theorem 3.1, there exists a unique graded Hopf pairing
  $\phi: S_{H_0}(H_1) \times S_{B_0}(B_1) \longrightarrow K$ extending $\phi_0$ and
  $\phi_1$.\par \vskip 3pt
    We claim that the following diagram commutes
  \[
\begin{CD}
H\times B @>  \psi>> K    \\
@V j_H \times j_B VV   @V {\rm Id}_K VV      \\
S_{H_0}(H_1) \times S_{B_0}{B_1} @> \phi >>K
\end{CD}\]
where the maps $j_H$ and $j_B$ are explained in Corollary 2.4.
\par

To see this, set $\psi'=\phi\circ (j_H \times j_B)$. Thus both
$\psi$ and $\psi'$ are graded Hopf pairings. Note that
\begin{align*}
\psi|_{H_0 \times B_0}= \psi'|_{H_0 \times B_0}\quad \mbox{and}
\quad \psi|_{H_1 \times B_1}= \psi'|_{H_1 \times B_1}.
\end{align*}
Since $H$ is generated by $H_0$ and $H_1$, it follows from the
fourth identity in the definition of Hopf pairing (see 3.1) that
$\psi=\psi'$. This shows the claim. \par

 By Corollary 2.4, the maps $j_H$ and $j_B$ are epimorphisms.
The fact that both $\phi$ and $\psi$ are two-sided non-degenerated
immediately implies that $j_H: H \simeq S_{H_0}(H_1)$ and $j_B: B
\simeq S_{B_0}(B_1)$. This completes the proof.\hfill
$\blacksquare$ \vskip15pt

\subsection{Self-dual couples} We end our paper with a special case of Theorem 3.1, which is of
independent interest. \par \vskip 5pt

Recall that a Hopf algebra $H$ is said to be self-dual, if there
exists a two-sided non-degenerated Hopf pairing $\phi: H \times H
\longrightarrow K$. Similarly, a graded Hopf algebra $H=\oplus_{n
\geq 0} H_n$ is said to be \textit{graded self-dual}, if the Hopf
pairing $\phi$ is graded. \par \vskip 5pt

A couple (H, M), where $H$ is a Hopf algebra and $M$ an $H$-Hopf
bimodule, is said to be \textit{self-dual}, if there exists a
pairing

  \begin{align*}
  (\phi_0, \phi_1): (H, M)\times (H, M) \longrightarrow K
  \end{align*}
such that both  $\phi_0$ and $\phi_1$ are two-sided
non-degenerated. Note that in this case, the $H$-Hopf bimodule $M$
is exactly the self-dual Hopf bimodule in [GM] and [HLY]. \par
\vskip 10pt

 The
following result is a direct consequence of Theorem 3.1.\par
\vskip 5pt

\begin{cor}
Let $(H, M)$ be a couple as above. Then the quantum symmetric
algebra $S_H(M)$ is graded self-dual if and only if the couple
$(H, M)$ is self-dual.
\end{cor}
\vskip30pt

\bibliography{}

\begin{thebibliography}{99999}
%\bibitem[ARS]{ARS} M. Auslander, I. Reiten, and S.O. Smal$\phi$,
%Representation Theory of Artin Algebras. Cambridge Studies in Adv.
%Math. 36, Cambridge Univ. Press,  1995.

\bibitem[CHYZ]{CHYZ} X. W. Chen, H. L. Huang, Y. Ye, and P. Zhang,
Monomial Hopf algebras, J. Algebra 275(2004), 212-232.

\bibitem[CHZ]{CHZ} X. W. Chen, H. L. Huang, and P. Zhang, Dual Gabriel theorem with applications,
Sci. in China, Ser. A Math. 49(1)(2006), 9-26.

\bibitem[CZ]{CZ} X. W. Chen and P. Zhang, Comodules of $U_q(sl_2)$ and modules of $SL_q(2)$
 via quiver methods, J. Pure Appl. Algebra, to appear.

\bibitem[CM]{CM} W. Chin and S. Montgomery, Basic coalgebras,
In: Modular interfaces (Reverside, CA, 1995), 41-47, AMS/IP Stud.
Adv. Math. 4, Amer. Math. Soc., Providence, RI, 1997.

\bibitem[C]{C}C. Cibils, A quiver quantum group, Comm. Math.
Phys. 157 (1993), 459-477.

\bibitem[CR1]{CR1} C. Cibils and M. Rosso, Algebres des chemins
quantique, Adv. Math. 125 (1997), 171-199.

\bibitem[CR2]{CR2} C. Cibils and
M. Rosso, Hopf quivers, J. Algebra 254(2002), 241-251.


%\bibitem[D]{D}Y. Doi, Homological coalgebras,
%J. Math. Soc. Japan 33(1)(1981), 31-50.


%\bibitem[DK]{DK}  Yu. A. Drozd, and V. V. Kirichenko, Finite Dimensional
%Algebras. Springer-Verlag, Berlin, Heidelberg, New York, Tokyo
%1993.

\bibitem[GM]{GM} E. L. Green and  E. N. Marcos, Self-dual Hopf
algebras, Comm. Algebra 28 (6) (2000), 2735-2744.

\bibitem[GS]{GS}E. L. Green and $\O$. Solberg,
Basic Hopf algebras and  quantum groups, Math. Z. 229(1998), 45-76.


\bibitem[HLY]{HLY}H. L. Huang,  L. B. Li, and Y. Ye, Self-dual Hopf
quivers, Comm. Algebra 33 (12) (2005), 4505-4514

\bibitem[K]{K} C. Kassel, Quantum Groups. Graduate Texts in
Math. 155, Springer-Verlag, New York, 1995.

\bibitem[M]{M}S.
Montgomery, Hopf Algebras and Their Actions on Rings. CBMS
Regional Conf. Series in Math. 82, Amer. Math. Soc., Providence,
RI, 1993.

\bibitem[N]{N} W. Nichols, Bialgebra of type I, Comm. Algebra 15
(1978), 1521-1552.



\bibitem[OZ]{OZ}F. van Oystaeyen  and P. Zhang, Quiver Hopf algebras,
J. Algebra 280(2004), 577-589.
%\bibitem[Rin]{Rin} C. M. Ringel, Tame Algebras and Integral Quadratic Forms. Lecture Notes in Math. 1099,
%Springer-Verlag, 1984.

\bibitem[Ro]{Ro} M. Rosso,  Quantum groups and quantum shuffles,
Invent. Math. 133 (1998), 339-416.

\bibitem[Sw]{Sw}M. E. Sweedler, Hopf Algebras, Benjamin, New York, 1969.
\end{thebibliography}

\end{document}